\documentclass[11pt,reqno]{amsart}
\usepackage{indentfirst, amssymb,amsmath,amsthm, mathrsfs,cite,graphicx,float}
\newtheorem{theo}{Theorem}[section]
\newtheorem{lem}{Lemma}[section]
\newtheorem{cor}{Corollary}[section]

\newtheorem{exam}{Example}[section]

\newtheorem{defi}{Definition}[section]

\newtheorem{open problem}{Open problem}[section]
\newcommand{\pa}{\partial}

\newcommand{\C}{\Bbb{C}}
\newcommand{\D}{\Bbb{D}}
\newcommand{\be}{\begin{equation}}
\newcommand{\ee}{\end{equation}}
\newcommand{\bs}{\begin{small}}
\newcommand{\es}{\end{small}}
\newcommand{\beas}{\begin{eqnarray*}}
\newcommand{\eeas}{\end{eqnarray*}}
\newcommand{\bea}{\begin{eqnarray}}
\newcommand{\eea}{\end{eqnarray}}
\renewcommand{\epsilon}{\varepsilon}
\numberwithin{equation}{section}

\begin{document}
\title[geometric investigation]{A geometric investigation of a certain subclass of univalent functions}
\author[R. Biswas and R. Mandal]{Raju Biswas and Rajib Mandal}
\date{}
\address{Raju Biswas, Department of Mathematics, Raiganj University, Raiganj, West Bengal-733134, India.}
\email{rajubiswasjanu02@gmail.com}
\address{Rajib Mandal, Department of Mathematics, Raiganj University, Raiganj, West Bengal-733134, India.}
\email{rajibmathresearch@gmail.com}

\maketitle
\let\thefootnote\relax

\begin{abstract} Let $\mathcal{H}$ be the space of all functions that are analytic in $\mathbb{D}$. Let $\mathcal{A}$ denote the family of all functions $f\in\mathcal{H}$ and normalized by the conditions $f(0)=0=f'(0)-1$. 
Obradovi\'{c} and Ponnusamy have introduced the class $\mathcal{M}(\lambda)$ such that the functions in $\mathcal{M}(\lambda)$ are univalent in $\mathbb{D}$ whenever $0<\lambda\leq 1$.
In this paper, we address a radius property of the class 
$\mathcal{M}(\lambda)$ and a number of associated results pertaining to $\mathcal{M}$.
The main objective of this paper is to examine the largest disks with sharp radius for which the functions $F$ defined by the relations $g(z)h(z)/z$, $z^2/g(z)$, and $z^2/\int_0^z (t/g(t))dt$ belong to 
the class $\mathcal{M}$, where $g$ and $h$ belong to some suitable subclasses of $\mathcal{S}$, the class of univalent functions from $\mathcal{A}$. 
In the final analysis, we obtain the sharp Bohr radius, Bohr-Rogosinski radius and improved Bohr radius for a certain subclass of starlike functions.
\end{abstract}
\section{ Introduction and preliminaries}
Let $\mathcal{H}$ be the class of all analytic functions in $\mathbb{D}$ and $\mathcal{A}$ denote the family of all functions $f\in\mathcal{H}$ and normalized by the 
conditions $f(0)=0=f'(0)-1$. Let $\mathcal{S}=\{f\in\mathcal{A}: f \;\text{is univalent in}\;\mathbb{D}\}$. 
The primary focus of this paper is on functions $f\in\mathcal{A}$ with $f(z)\not=0$ for $z\in\Bbb{D}$ of the form 
\bea\label{y1} \frac{z}{f(z)}=1+\sum_{n=1}^\infty b_n z^n, \quad z\in\Bbb{D}. \eea
Obviously, if $f\in\mathcal{S}$, then $z/f(z)\not=0$ in $\Bbb{D}$, and thus it can be expressed as a Taylor series expansion of the form (\ref{y1}).  We say that a set 
$E\subseteq \C$ is said to be starlike with respect to a point $w_0\in E$ if the linear segment joining $w_0$ to every other point $w\in E$ lies entirely in $E$, {\it i.e.,} 
$t w_0+(1-t)w\in E$, where $w\in E$ and $t\in[0,1]$. A function $f\in\mathcal{S}$ is called starlike with respect to origin if $f(\D)$ is starlike with respect to origin. We denote this 
set of functions by $\mathcal{S}^*$. Analytically, if $f\in \mathcal{S}^*$ if, and only if, $\text{Re}\left(zf'(z)/f(z)\right)>0$ in $\D$ (see \cite[Theorem 2.10, P. 41]{2}).\\[2mm]
\indent Let $\mathcal{U}(\lambda)$ denote the class of all functions $f\in\mathcal{A}$ in $\mathbb{D}$ satisfying the condition 
\beas \left|f'(z)\left(\frac{z}{f(z)}\right)^2-1\right|\leq \lambda\quad \text{for}\quad z\in\mathbb{D}\quad\text{with}\quad\lambda>0.\eeas
We denote the class $\mathcal{U}(1)$ by $\mathcal{U}$. Moreover $\mathcal{U}(\lambda)\subset\mathcal{U}$ for $0<\lambda\leq1$ (see \cite{1,1a,14}). The functions in 
the class $\mathcal{U}(\lambda)$ are univalent in $\mathbb{D}$ whenever $0<\lambda\leq 1$. In a more general sense, $\mathcal{U}(\lambda)\subsetneq\mathcal{S}$ for 
$ 0<\lambda\leq1$ (see \cite{5}). 
The properties of the class $\mathcal{U}$ have been the subject of extensive study in \cite{3,5,7,10,12,19,6,9,11,8,17,18}. It is well established that there are only nine functions in $\mathcal{S}$ with integral coefficients in the power series expansions of $f\in\mathcal{S}$ (see \cite{20}). If $\mathcal{S}_{\mathbb{Z}}=\{f\in\mathcal{S}: a_n\in\mathbb{Z}\}$, then
\beas \mathcal{S}_{\mathbb{Z}}=\left\{z,\frac{z}{(1\pm z)^2}, \frac{z}{1\pm z}, \frac{z}{1\pm z^2}, \frac{z}{1\pm z+z^2}\right\}.\eeas
In 2011, Obradovi\'{c} and Ponnusamy\cite{13} introduced the following important class:
\beas \mathcal{M}(\lambda)=\left\{f\in\mathcal{A}: \left|z^2\left(\frac{z}{f(z)}\right)''+f'(z)\left(\frac{z}{f(z)}\right)^2-1\right|\leq \lambda\quad\text{for}\quad z\in\mathbb{D}, \lambda>0\right\}.\eeas 
We denote the class $\mathcal{M}(1)$ by $\mathcal{M}$ and the class $\mathcal{M}_2$ by $\left\{f\in\mathcal{M}: f''(0)=0\right\}$. 
Obradovi\'{c} and Ponnusamy \cite{13} have shown that for $0<\lambda\leq1$, the strict inclusion relation 
$\mathcal{M}(\lambda)\subsetneq \mathcal{U}(\lambda)\cap\mathcal{P}(\lambda)\subsetneq \mathcal{S}$ holds, where 
$\mathcal{P}(\lambda)=\left\{f\in\mathcal{A}: \left\vert\left(z/f(z)\right)''\right\vert\leq 2\lambda, \;z\in\mathbb{D}\right\}$. Functions in $\mathcal{M}$ are also known to be 
univalent in $\mathbb{D}$. The class $\mathcal{M}$ is closely related to the class $\mathcal{U}$ in the sense of the strict inclusion 
$\mathcal{S}_{\mathbb{Z}}\subsetneq\mathcal{M}\subsetneq\mathcal{P}\subsetneq\mathcal{U}\subsetneq\mathcal{S}$ (see \cite{1,15,13,14}). Consequently, it has been 
verified that each function in $\mathcal{S}_{\mathbb{Z}}$ belongs to $\mathcal{U}\cap\mathcal{M}$ (see \cite{13}). It is well known that the class $\mathcal{S}$ is preserved under a number of elementary transformations (see \cite[P. 27]{2}). But, it is easy to check that the class $\mathcal{M}(\lambda)$ is not preserved under rotation, conjugation and dilation, while $\mathcal{M}(\lambda)$ is preserved under omitted-value transformation. Another concept that has been widely discussed is the 
Hankel determinant of the logarithmic coefficients of univalent functions. See \cite{V1,103} and the references therein for some recent results on this topic.
We consider the following subclasses of $\mathcal{A}$:
\beas \mathcal{R}(1/2)&=&\left\{f\in\mathcal{A}: \text{Re}\left(\frac{f(z)}{z}\right)>\frac{1}{2}, \quad z\in\Bbb{D}\right\},\\[2mm]
\mathcal{C}(-1/2)&=&\left\{f\in\mathcal{A}: \text{Re}\left(1+\frac{zf''(z)}{f'(z)}\right)>-\frac{1}{2},\quad  z\in\Bbb{D}\right\}\quad\text{and}\\[2mm]
\mathcal{G}&=&\left\{f\in\mathcal{A}: \text{Re}\left(1+\frac{zf''(z)}{f'(z)}\right)<\frac{3}{2},\quad  z\in\Bbb{D}\right\}.
\eeas
It is known that $\mathcal{K}\subset \mathcal{R}(1/2)$ (see \cite[P. 67]{r3}), where $\mathcal{K}$ denotes the class of all functions $f\in\mathcal{A}$ that are convex in 
$\Bbb{D}$. The classes $\mathcal{C}(-1/2)$ and $\mathcal{G}$ have been extensively studied in \cite{r4,r5}. Further, functions belonging to the class $\mathcal{C}(-1/2)$ are 
known to be close-to-convex in $\Bbb{D}$. Ozaki \cite{r6} introduced the class $\mathcal{G}$ and proved that the functions in $\mathcal{G}$ are univalent in $\Bbb{D}$. However, 
the functions belonging to the class $\mathcal{G}$ are starlike in $\Bbb{D}$ (see \cite[Theorem 1]{r7} and \cite{r5}).\\[2mm]
\indent For the purposes of this study, the following definition is necessary.
\begin{defi}  The polylogarithm $Li_k(z)$ is a special function of order $k$ and argument $z$ defined in the complex plane over the unit disk $\Bbb{D}$ as
\beas \text{Li}_k(z)=\sum_{n=1}^\infty\frac{z^n}{n^k}=z+\frac{z^2}{2^k}+\frac{z^3}{3^k}+\cdots\eeas
and it is also known as Jonqui\`ere's function. The special cases $k=2$ and $k=3$ are called the dilogarithm and trilogarithm respectively. \end{defi}
The remaining sections are organized as follows: In section $2$, we present a number of lemmas and a proof of one of these lemmas, which are essential for the proof
of our main theorems. 
In section $3$, we investigate a radius property of the class $\mathcal{M}(\lambda)$ and a number of associated results pertaining to $\mathcal{M}$. In section $4$, we examine 
the largest disks for which the functions defined by the quotient of analytic functions belongs to $\mathcal{M}$. Sections $5$ and 
$6$ comprise the introductory sections of a certain subclass of starlike functions and the Bohr phenomenon regarding the class of bounded analytic functions, respectively. 
In section $7$, we obtain the sharp Bohr radius, Bohr-Rogosinski radius and improved Bohr radius for a certain subclass of starlike functions.
\section{Necessary lemmas} 
\begin{figure}[H]
\begin{minipage}[c]{0.45\linewidth}
\centering
\includegraphics[scale=0.5]{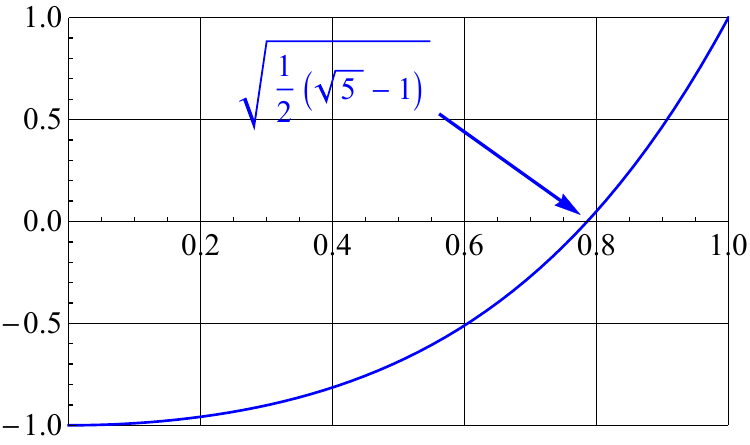}
\caption{The graph of the polynomial $r^4+r^2-1$ in $(0,1)$ }
\label{fig2}
\end{minipage}\hspace{0.5cm}
\begin{minipage}[c]{0.45\linewidth}
\centering
\includegraphics[scale=0.5]{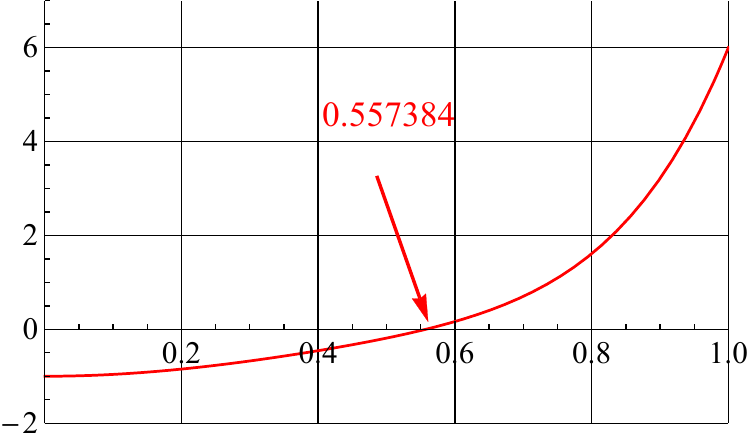}
\caption{The graph of the polynomial $8r^6- 5r^4 +4r^2-1$ in $(0,1)$}
\label{fig1}
\end{minipage}
\end{figure}

For this paper, the following lemmas will be used to prove the main results.
\begin{lem}\label{lem}\cite[Theorem 4]{13} Every $f\in\mathcal{M}(\lambda)$ has the representation
\beas \frac{z}{f(z)}=1-\frac{f''(0)}{2}z+\lambda \int_0^1\frac{\omega(tz)}{t^2}\log(1/t)dt,\eeas
for some $\omega: \Bbb{D}\to\Bbb{D}$ with $\omega(0)=\omega'(0)=0$.
\end{lem}
\begin{lem}\cite[Area Theorem, Theorem 11 of Vol. 2]{4}\label{lem0} Let $\mu>0$ and $f\in\mathcal{S}$ be in the form $\left(z/f(z)\right)^\mu=1+\sum_{n=1}^\infty b_nz^n$. Then, we have $\sum_{n=1}^\infty (n-\mu)|b_n|^2\leq \mu$.
\end{lem}
\begin{lem}\label{lem1}\cite[Theorem 2, p. 171]{13} Let $\phi(z)=1+\sum_{n=1}^\infty b_nz^n$ be a non-vanishing analytic function in $\mathbb{D}$ that satisfying the coefficient bound $\sum_{n=2}^\infty (n-1)^2|b_n|\leq \lambda$.
Then the function $f$ defined by $f(z)=z/\phi(z)$ is in $\mathcal{M}(\lambda)$.\end{lem}
\begin{lem}\label{lem2}\cite[Corollary 1, p. 172]{13} Let $f\in\mathcal{M}$ be of the form $z/f(z)=1+\sum_{n=1}^\infty b_nz^n$. Then, we have $\sum_{n=2}^\infty (n-1)^4|b_n|^2\leq 1$.\end{lem}
It is evident that the class $\mathcal{M}$ is not preserved under dilation. This motivates the determination of the largest disk in which the class $\mathcal{M}$ is preserved under dilation.
\begin{lem}\label{T3} If $f\in\mathcal{M}$ and $g(z)=f(rz)/r$, then $g\in\mathcal{M}$ for $0<r\leq r_0$, where $r_0=\sqrt{(\sqrt{5}-1)/2}\approx 0.786151$ is the unique root of the equation $r^4+r^2-1=0$ in $(0,1)$.\end{lem}
\begin{proof}
Let $f\in\mathcal{M}$. Then, $f$ is univalent and $f(z)/z$ can be expressed as
\bea\label{eq8} \frac{z}{f(z)}=1+\sum_{n=1}^\infty b_nz^n\quad\text{for}\quad z\in\mathbb{D}.\eea
Let us consider the function $g(z)=f(rz)/r$ for $0<r<1$. From (\ref{eq8}), we have 
\beas \frac{z}{f(rz)/r}=1+\sum_{n=1}^\infty b_nr^nz^n.\eeas  
We need to show that $f(rz)/r\in\mathcal{M}$ for $0<r\leq r_0$. 
In view of \textrm{Lemma \ref{lem1}}, it is sufficient to show that $\sum_{n=2}^\infty (n-1)^2\left|b_nr^n\right|\leq1$ for $0<r\leq r_0$.
Using the Cauchy-Schwarz inequality, we have 
\beas\sum_{n=2}^\infty (n-1)^2\left|b_nr^n\right|\leq \left(\sum_{n=2}^\infty (n-1)^4\left|b_n\right|^2\right)^{1/2}\left(\sum_{n=2}^\infty r^{2n}\right)^{1/2}.\eeas
Since $f\in\mathcal{M}$, it follows from \textrm{Lemma \ref{lem2}} and relation (\ref{eq8}) that 
$\sum_{n=2}^\infty (n-1)^4\left|b_n\right|^2\leq 1$.
Therefore, we have
\beas\sum_{n=2}^\infty (n-1)^2\left|b_nr^n\right|\leq\left(\sum_{n=2}^\infty r^{2n}\right)^{1/2}=\left(\frac{r^4}{\left(1-r^2\right)}\right)^{1/2}\leq 1\eeas
for $0<r\leq r_0$, where $r_0\left(=\sqrt{(\sqrt{5}-1)/2}\right)$ is the unique root of the equation $r^4+r^2-1=0$ in $(0,1)$, as illustrated in Figure \ref{fig2}. This completes the proof.
\end{proof}
\section{Results on Radius properties}
\noindent Suppose $\mathcal{K}_1,\mathcal{K}_2\subset\mathcal{A}$. If for every $f\in\mathcal{K}_1$, there exists a largest number $r_0$ such that $f(rz)/r\in\mathcal{K}_2$ for 
$r\leq r_0$, then $r_0$ is called the $\mathcal{K}_2-$radius in $\mathcal{K}_1$. A substantial body of research exists in the field of univalent function theory, with numerous 
studies exploring this topic (see \cite{9,11,13,19}). In the following, we establish the $\mathcal{M}(\lambda)$-radius in $\mathcal{S}$.
\begin{theo}\label{Theo} Let $f\in\mathcal{S}$ and $\lambda\in(0,1]$. Then the function $f(rz)/r$ belongs to $\mathcal{M}(\lambda)$ for $0<r\leq r_0$, where $r_0\in(0,1)$ is the 
smallest root the smallest root of the equation
\beas r^4(r^4+4r^2+1)-\lambda^2\left(1-r^2\right)^4=0.\eeas\end{theo}
\begin{proof} As $f\in\mathcal{S}$, so $z/f(z)\not=0$ in $\mathbb{D}$. Then, $z/f(z)$ can be written as
\bea\label{eqr8} \frac{z}{f(z)}=1+\sum_{n=1}^\infty b_nz^n\;\;\text{for}\;z\in\mathbb{D}.\eea
In view of \textrm{Lemma \ref{lem0}}, we have $\sum_{n=2}^\infty (n-1)|b_n|^2\leq 1$.
Let $g(z)=f(rz)/r$ for $0<r<1$. From (\ref{eqr8}), we have 
\beas \frac{z}{f(rz)/r}=1+\sum_{n=1}^\infty b_nr^nz^n.\eeas  
We need to show that $f(rz)/r\in\mathcal{M}(\lambda)$ for $0<r\leq r_0$. 
In view of \textrm{Lemma \ref{lem1}}, it is sufficient to show that $\sum_{n=2}^\infty (n-1)^2\left|b_nr^n\right|\leq\lambda$ for $0<r\leq r_0$.
In view of Cauchy-Schwarz inequality, we have 
\beas\sum_{n=2}^\infty (n-1)^2\left|b_nr^n\right|&\leq& \left(\sum_{n=2}^\infty (n-1)\left|b_n\right|^2\right)^{1/2}\left(\sum_{n=2}^\infty (n-1)^{3}r^{2n}\right)^{1/2}\\[2mm]
&\leq&\left(\sum_{n=2}^\infty (n-1)^{3}r^{2n}\right)^{1/2}=\frac{r^2\left(r^4+4r^2+1\right)^{1/2}}{\left(1-r^2\right)^2}\leq \lambda\eeas
for $0<r\leq r_0$, where $r_0\in(0,1)$ is the smallest root of the equation 
\beas r^4(r^4+4r^2+1)-\lambda^2\left(1-r^2\right)^4=0.\eeas
This completes the proof.
\end{proof}
\begin{cor} Suppose $f\in\mathcal{S}$. Then, the function $f(rz)/r$ belongs to $\mathcal{M}$ for $0<r\leq r_0$, where $r_0\thickapprox 0.557384$ is the unique root of the equation $8r^6- 5r^4 +4r^2-1=0$ in $(0,1)$, as shown in Figure \ref{fig1}.
\end{cor}
\begin{proof} The result immediately follows from \textrm{Theorem \ref{Theo}}.\end{proof}
\section{Quotient of analytic functions}
\begin{figure}[H]
\begin{minipage}[c]{0.45\linewidth}
\centering
\includegraphics[scale=0.5]{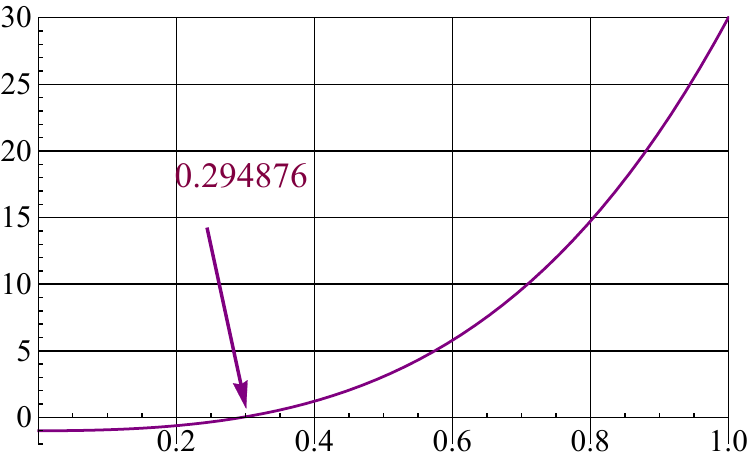}
\caption{The graph of the polynomial $9r^4+16r^3+6r^2-1$ in $(0,1)$}
\label{fig3}
\end{minipage}\hspace{1cm}
\begin{minipage}[c]{0.45\linewidth}
\centering
\includegraphics[scale=0.5]{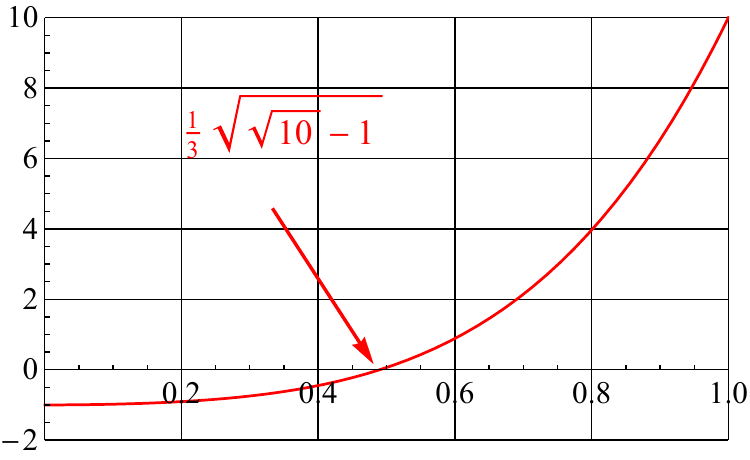}
\caption{The graph of the polynomial $9r^4+2r^2-1$}
\label{fig4}
\end{minipage}
\end{figure}

\begin{figure}[H]
\includegraphics[scale=0.5]{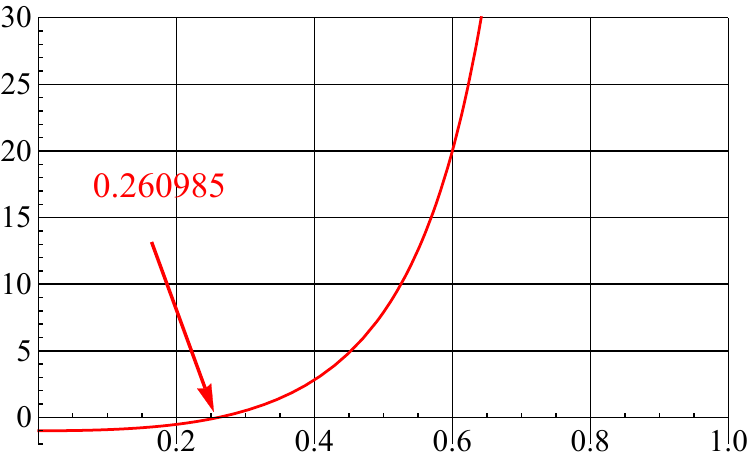}
\caption{The graph $A(r)$ in $(0,1)$}
\label{fig5}
\end{figure}

\begin{figure}[H]
\begin{minipage}[c]{0.4\linewidth}
\centering
\includegraphics[scale=0.55]{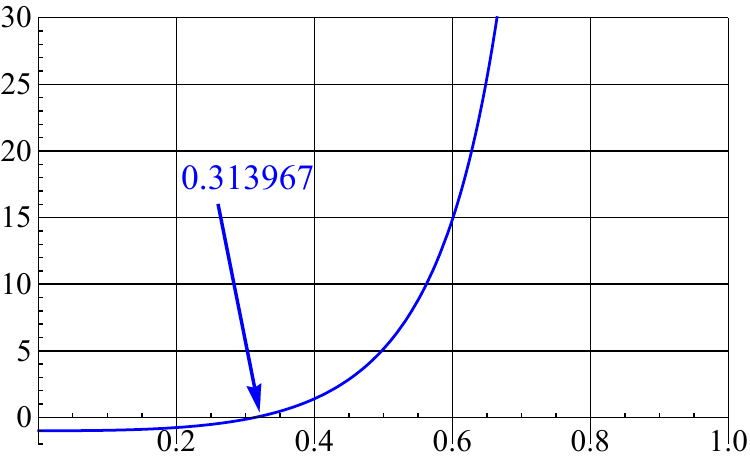}
\caption{The graph of $B(r)$ in $(0,1)$}
\label{fig6}
\end{minipage}\hspace{1cm}
\begin{minipage}[c]{0.4\linewidth}
\centering
\includegraphics[scale=0.55]{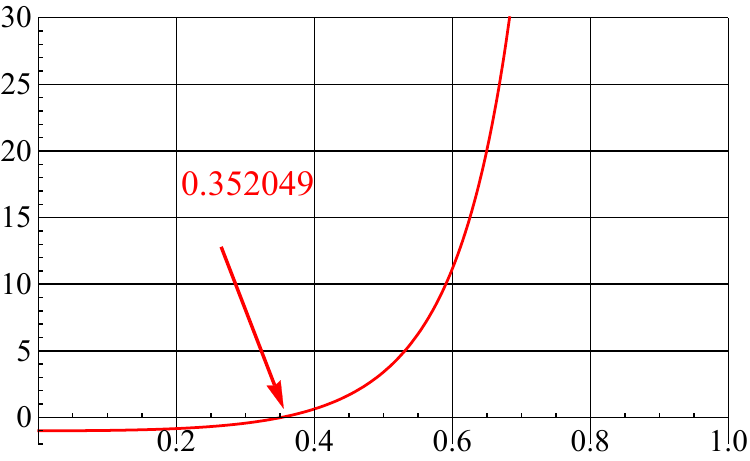}
\caption{The graph of $C(r)$ in $(0,1)$}
\label{fig7}
\end{minipage}
\end{figure}

\begin{figure}[H]
\begin{minipage}[c]{0.45\linewidth}
\centering
\includegraphics[scale=0.5]{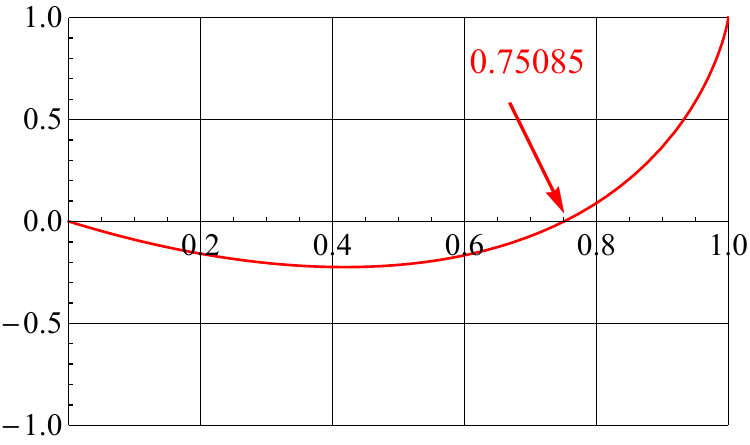}
\caption{The graph of $3r-2r^2+(r^2- 5 r+4)\log(1-r)$ in $(0,1)$}
\label{fig8}
\end{minipage}\hspace{0.5cm}
\begin{minipage}[c]{0.45\linewidth}
\centering
\includegraphics[scale=0.5]{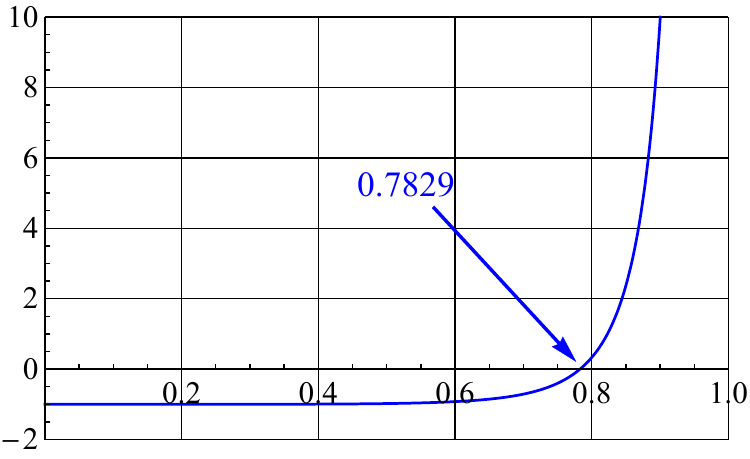}
\caption{The graph of $A_3(r)$ in $(0,1)$}
\label{fig9}
\end{minipage}
\end{figure}

\noindent Motivated by the study of Obradovi\'c and Ponnusamy \cite{8}, in this paper, we consider the following quotient of analytic functions:\\[2mm]
For $g\in\mathcal{G}_1\subseteq\mathcal{S}$ and $h\in\mathcal{G}_2\subseteq\mathcal{S}$, we consider the quotient function $F$ defined by
\bea\label{rm} F(z)=\frac{g(z)h(z)}{z}\quad\text{for}\quad z\in\mathbb{D}.\eea
For suitable choices of $\mathcal{G}_1$ and $\mathcal{G}_2$, our objective is to determine the largest disk within which $F\in\mathcal{M}$.
\begin{theo}\label{T1} Suppose that $g,h\in\mathcal{M}$. Then, the quotient function $F$ defined by (\ref{rm}) belongs to $\mathcal{M}$ in the disk $|z|\leq r_0$, where $r_0(\approx 0.294876)$ is the smallest root of the equation $9r^4+16r^3+6r^2-1=0$ in $(0,1)$. The result is sharp.\end{theo}
\begin{proof} Since $g\in\mathcal{M}\subsetneq\mathcal{U}$, thus, we have
\bea&&\label{eq1} -z\left(\frac{z}{g(z)}\right)'+\frac{z}{g(z)}-1=\omega(z),\eea 
where $\omega:\mathbb{D}\to\mathbb{D}$ is analytic and $\omega(0)=\omega'(0)=0$. From (\ref{eq1}), we have
\bea\label{eq2}\left|g'(z)\left(\frac{z}{g(z)}\right)^2-1\right|\leq|z|^2.\eea
Furthermore, as $g\in\mathcal{M}$, we have
\beas z^2\left(\frac{z}{g(z)}\right)''+g'(z)\left(\frac{z}{g(z)}\right)^2-1=\omega_1(z),\eeas
where $\omega_1:\mathbb{D}\to\mathbb{D}$ is analytic with $\omega_1(0)=\omega_1'(0)=0$. In view of the classical Schwarz lemma, we have $|\omega_1(z)|\leq |z|^2$.  Using \textrm{Lemma \ref{lem}}, we have
\beas \frac{z}{g(z)}=1-b_1z+\int_{0}^1\frac{\omega_1(tz)}{t^2}\log\left(\frac{1}{t}\right)dt,\quad\text{where}\quad b_1=\frac{g''(0)}{2!}.\eeas  
Thus, we have 
\bea\label{eq3} \left|z^2\left(\frac{z}{g(z)}\right)''+g'(z)\left(\frac{z}{g(z)}\right)^2-1\right|\leq|z|^2\quad\text{and}\quad\left|\frac{z}{g(z)}-1\right|\leq |b_1||z|+|z|^2.\eea
A similar conclusion holds for $h\in\mathcal{M}$. Since each functions in $\mathcal{M}$ are univalent, the Bieberbach's theorem of the univalent function gives that $|g''(0)/2|=|b_1|\leq2$ (see \cite[P. 30]{2} and \cite[P. 4]{r3}).
From the definition of $F$, we  have $F(z)/z\not=0$ in $\mathbb{D}$, since $g,h\in\mathcal{M}$.
Through straightforward mathematical calculations, we obtain that
\beas F'(z)\left(\frac{z}{F(z)}\right)^2-1&=&\left(g'(z)\left(\frac{z}{g(z)}\right)^2-1\right)\frac{z}{h(z)}+\left(h'(z)\left(\frac{z}{h(z)}\right)^2-1\right)\frac{z}{g(z)}\\[2mm]
&&-\left(\frac{z}{g(z)}-1\right)\left(\frac{z}{h(z)}-1\right)\quad\text{and}\\[2mm]
z^2\left(\frac{z}{F(z)}\right)''&=&2\left(\frac{z}{F(z)}\right)^3\left(F'(z)\right)^2-\left(\frac{z}{F(z)}\right)^2zF''(z)-2\left(\frac{z}{F(z)}\right)^2F'(z)\\
&=&\frac{z}{h(z)}\left(z^2\left(\frac{z}{g(z)}\right)''+g'(z)\left(\frac{z}{g(z)}\right)^2-1\right)\\
&&+\frac{z}{g(z)}\left(z^2\left(\frac{z}{h(z)}\right)''+h'(z)\left(\frac{z}{h(z)}\right)^2-1\right)\\
&&+2\left(g'(z)\left(\frac{z}{g(z)}\right)^2-1\right)\left(h'(z)\left(\frac{z}{h(z)}\right)^2-1\right)\\
&&-2\left(g'(z)\left(\frac{z}{g(z)}\right)^2-1\right)\left(\frac{z}{h(z)}-1\right)\\
&&-2\left(h'(z)\left(\frac{z}{h(z)}\right)^2-1\right)\left(\frac{z}{g(z)}-1\right)\\
&&-\frac{z}{h(z)}\left(g'(z)\left(\frac{z}{g(z)}\right)^2-1\right)-\frac{z}{g(z)}\left(h'(z)\left(\frac{z}{h(z)}\right)^2-1\right)\\
&&+2\left(\frac{z}{g(z)}-1\right)\left(\frac{z}{h(z)}-1\right).\eeas
Our objective is to find the largest disk $|z|\leq r$ for which the condition
\beas\left|z^2\left(\frac{z}{F(z)}\right)''+F'(z)\left(\frac{z}{F(z)}\right)^2-1\right|\leq 1\quad \text{holds}.\eeas 
It is evident that the condition $\left|z^2\left(z/F(z)\right)''+F'(z)\left(z/F(z)\right)^2-1\right|\leq 1$ holds in the disk $|z|\leq r$ if the inequality 
\bs\bea\label{eq4} &&\left|\frac{z}{h(z)}\right|\left|z^2\left(\frac{z}{g(z)}\right)''+g'(z)\left(\frac{z}{g(z)}\right)^2-1\right|+\left|\frac{z}{g(z)}\right|\left|z^2\left(\frac{z}{h(z)}\right)''+h'(z)\left(\frac{z}{h(z)}\right)^2-1\right|\nonumber\\[2mm]
&&+2\left|g'(z)\left(\frac{z}{g(z)}\right)^2-1\right|\left|h'(z)\left(\frac{z}{h(z)}\right)^2-1\right|+2\left|g'(z)\left(\frac{z}{g(z)}\right)^2-1\right|\left|\frac{z}{h(z)}-1\right|\nonumber\\[2mm]
&&+2\left|h'(z)\left(\frac{z}{h(z)}\right)^2-1\right|\left|\frac{z}{g(z)}-1\right|+\left|\frac{z}{g(z)}-1\right|\left|\frac{z}{h(z)}-1\right|\leq 1\eea\es
holds in the disk $|z|\leq r$. 
By employing the conditions (\ref{eq2}) and (\ref{eq3}), we are led to the conclusion that the inequality (\ref{eq4}) holds if
\beas
&&9|z|^4+16|z|^3+6|z|^2-1\leq 0.\eeas
Therefore, the function $F$ belongs to the class $\mathcal{M}$ for the disk $|z|\leq r_0$, where $r_0(\approx 0.294876)$ is the unique root of the equation $ F(r)=9r^4+16r^3+6r^2-1=0$ in $(0,1)$ and it is illustrated in Figure \ref{fig3}. It is evident that the function $F(r)$ is monotonically increasing with respect to $r$, with $F(0) = -1 < 0$ and $F(1)=30$. The Intermediate Value Theorem implies the existence of a unique $r_0\in(0, 1)$ such that $F(r_0)=0$.\\[2mm]
\indent To prove the sharpness of the radius $r_0$, we consider $g(z)=h(z)=z/(1-z)^2$. Then, $g, h\in\mathcal{M}$ and the corresponding quotient function $F$ gives that
\beas\left|z^2\left(\frac{z}{F(z)}\right)''+F'(z)\left(\frac{z}{F(z)}\right)^2-1\right|=|z|^2\left|9z^2-16z+6\right|. \eeas
Therefore, 
\beas \left|z^2\left(\frac{z}{F(z)}\right)''+F'(z)\left(\frac{z}{F(z)}\right)^2-1\right|_{z=-r}=\left(9r^4+16r^3+6r^2-1\right)+1>1\eeas
for $r>r_0$, which shows that the number $r_0$ is sharp. This completes the proof.
 \end{proof}
\begin{cor} Suppose that $g,h\in\mathcal{M}_2$. Then, the quotient function $F(z)=g(z)h(z)/z$ for $z\in\mathbb{D}$ belongs to $\mathcal{M}_2$ in the disk $|z|\leq r_0$, where $r_0(=\sqrt{\sqrt{10}-1}/3)$ is the unique root of the equation $9r^4+2r^2-1=0$ in $(0,1)$, as illustrated in Figure \ref{fig4}.\end{cor}
\begin{proof} 
The result is an immediate consequence of \textrm{Theorem \ref{T1}}, with the additional conditions $b_1=c_1=0$.\end{proof}
\begin{cor} Suppose that $g\in\mathcal{M}_2(\lambda)$ and $h\in\mathcal{M}_2(\lambda')$, where $\lambda,\lambda'\in(0,1]$. Then, the quotient function $F$ defined by (\ref{rm}) belongs to $\mathcal{M}_2(\mu)$ in the disk $|z|\leq r_1$, where 
\beas r_1=\sqrt{\frac{-(\lambda+\lambda')+\sqrt{(\lambda+\lambda')^2+12\mu(\lambda^2+\lambda\lambda'+\lambda'^2)}}{6(\lambda^2+\lambda\lambda'+\lambda'^2)}}.\eeas\end{cor}
\begin{proof} 
The result is an immediate consequence of \textrm{Theorem \ref{T1}}.\end{proof}
\begin{theo}\label{T4}
Let $g,h\in\mathcal{S}$. Then, the quotient function $F$ defined by (\ref{rm}) belongs to the class $\mathcal{M}$ in the disk $|z|\leq r_0$, where $r_0(\approx0.260985)$ is the 
smallest positive root of the equation
\beas &&6r^2+4(\sqrt{2}+4)r^3+\frac{2r^4\sqrt{-8r^{10}+31r^8-44r^6+27r^4}}{(1-r^2)^{2}}+\frac{r^4(4r^2-11r+9)}{(1-r)^3}\\
&&+4\left(\frac{r^8(15-20 r^2+15 r^4-4 r^6)}{(1-r^2)^4}-r^4\left(\log(1-r^2)+r^2\right)\right)^{1/2}-1=0\eeas
in $(0,1)$.
\end{theo}
\begin{proof}
Since $g,h\in\mathcal{S}$, so $z/g(z)\not=0$ and $z/h(z)\not=0$ in $\mathbb{D}$. Then, $z/g$ and $z/h$ can be expressed as
\beas \frac{z}{g(z)}=1+\sum_{n=1}^\infty b_nz^n\quad\text{and}\quad\frac{z}{h(z)}=1+\sum_{n=1}^\infty c_nz^n\quad\text{for}\quad z\in\mathbb{D},\eeas
where $b_1=-g''(0)/2$ and $c_1=-h''(0)/2$. In view of the Bieberbach's theorem, we have $|b_1|\leq 2$ and $|c_1|\leq 2$. Since $g, h\in\mathcal{S}$, it follows from \textrm{Lemma \ref{lem0}}, we have 
\bea\label{eq5} \sum_{n=2}^\infty (n-1)|b_n|^2\leq 1\quad\text{and}\quad\sum_{n=2}^\infty (n-1)|c_n|^2\leq 1.\eea
In particular,
\bea\label{eq6} \sum_{n=2}^\infty|b_n|^2\leq \sum_{n=2}^\infty (n-1)|b_n|^2\leq 1\quad\text{and}\quad\sum_{n=2}^\infty |c_n|^2\leq 1.\eea
Therefore
\bea\label{eq7} \frac{z}{F(z)}=\left(1+\sum_{n=1}^\infty b_nz^n\right)\left(1+\sum_{n=1}^\infty c_nz^n\right)=1+\sum_{n=1}^\infty B_nz^n,\eea
where $B_n=\sum\limits_{i=0}^n b_ic_{n-i}$ with $b_0=c_0=1$. From (\ref{eq5}), we have 
\beas |b_2|^2\leq 1\Rightarrow |b_2|\leq1\quad\text{and}\quad 2|b_3|^2\leq 1\Rightarrow |b_3|\leq \frac{1}{\sqrt{2}}.\eeas
 Similarly, we deduce that $|c_2|\leq1$ and $|c_3|\leq 1/\sqrt{2}$. It is evident that
 \beas |B_2|\leq 2+|b_1c_1|\quad\text{and}\quad |B_3|\leq \sqrt{2}+|b_1|+|c_1|.\eeas
Using a consequences of Cauchy-Schwarz inequality and (\ref{eq6}), we have  
\beas  |B_n|&\leq& |b_n|+|c_n|+|b_1||c_{n-1}|+|c_1||b_{n-1}|+\sum_{i=2}^{n-2}|b_i||c_{n-i}|\\[2mm]
&\leq&|b_n|+|c_n|+|b_1||c_{n-1}|+|c_1||b_{n-1}|+\left(\sum_{i=2}^{n-2}|b_i|^2\right)^{1/2}\left(\sum_{i=2}^{n-2}|c_i|^2\right)^{1/2}\\[2mm]
&\leq&|b_n|+|c_n|+|b_1||c_{n-1}|+|c_1||b_{n-1}|+1\quad\text{for}\quad n\geq 4.\eeas
We consider the function $G(z)=F(rz)/r$ for $0<r<1$. By (\ref{eq7}), we have 
\beas \frac{z}{G(z)}=1+\sum_{n=1}^\infty B_nr^nz^n.\eeas
Our objective is to find the largest disk $|z|\leq r$ for which $G\in\mathcal{M}$. In view of \textrm{Lemmas \ref{lem1}} and \ref{T3}, it is sufficient to show that
\beas S=\sum_{n=2}^\infty (n-1)^2|B_n|r^n=|B_2|r^2+4|B_3|r^3+S_1\leq 1\;\; \text{for}\;\;0<r\leq \sqrt{(\sqrt{5}-1)/2},\eeas
where 
\beas S_1=\sum_{n=4}^\infty (n-1)^2|B_n|r^n&\leq& \sum_{n=4}^\infty (n-1)^2\left(|b_n|+|c_n|+|b_1||c_{n-1}|+|c_1||b_{n-1}|+1\right)r^n\\[2mm]
&=&T_1+T_2+|b_1|T_3+|c_1|T_4+T_5\quad\text{with}\eeas
\beas && T_1=\sum_{n=4}^\infty (n-1)^2|b_n|r^n,\quad T_2=\sum_{n=4}^\infty (n-1)^2|c_n|r^n,\quad T_3=\sum_{n=4}^\infty (n-1)^2|c_{n-1}|r^n,\\[2mm]
&& T_4=\sum_{n=4}^\infty (n-1)^2|b_{n-1}|r^n\quad\text{and}\quad T_5=\sum_{n=4}^\infty (n-1)^2r^n=\frac{r^4(4r^2-11r+9)}{(1-r)^3}.\eeas
It is evident that
\beas S\leq (2+|b_1c_1|)r^2+4(\sqrt{2}+|b_1|+|c_1|)r^3+S_1.\eeas
Using (\ref{eq5}), (\ref{eq6}) and the Cauchy-Schwarz inequality, we have
\beas T_1= \sum_{n=4}^\infty (n-1)^2|b_n|r^n&\leq& \left(\sum_{n=4}^\infty (n-1)|b_n|^2\right)^{1/2}\left(\sum_{n=4}^\infty (n-1)^3r^{2n}\right)^{1/2}\\[2mm]
&\leq&\left(\sum_{n=4}^\infty (n-1)^3r^{2n}\right)^{1/2}\\[2mm]
&=&\frac{r^2\sqrt{-8r^{10}+31r^8-44r^6+27r^4}}{(1-r^2)^{2}}\quad\text{and}\eeas
\beas T_3=  \sum_{n=4}^\infty (n-1)^2|c_{n-1}|r^n&\leq& \left(\sum_{n=4}^\infty (n-2)|b_{n-1}|^2\right)^{1/2}\left(\sum_{n=4}^\infty \frac{(n-1)^4r^{2n}}{n-2}\right)^{1/2}\\[2mm]
&\leq&\left(\sum_{n=4}^\infty \frac{(n-1)^4r^{2n}}{n-2}\right)^{1/2}.\eeas
Now
\beas\sum_{n=4}^\infty \frac{(n-1)^4r^{n}}{n-2}&=&\sum_{n=4}^\infty (n-2)^3r^n+4\sum_{n=4}^\infty (n-2)^2r^n+6\sum_{n=4}^\infty (n-2)r^n+4\sum_{n=4}^\infty r^n\\[2mm]
&&+\sum_{n=4}^\infty \frac{r^n}{(n-2)}\\[2mm]
&=&\left(\frac{r^2\left(r^3+4r^2+r\right)}{(1-r)^4}-r^3\right)+4\left(\frac{r^2\left(r^2+r\right)}{(1-r)^3}-r^3\right)\\[2mm]
&&+6r^2\left(\frac{r}{(1-r)^2}-r\right)+\frac{4r^4}{1-r}-r^2(\log(1-r)+r)\\[2mm]
&=&\frac{r^4\left(15-20 r+15 r^2-4 r^3\right)}{(1-r)^4}-r^2(\log(1-r)+r).\eeas
Thus, we have
\beas T_3\leq \left(\frac{r^8\left(15-20 r^2+15 r^4-4 r^6\right)}{(1-r^2)^4}-r^4\left(\log(1-r^2)+r^2\right)\right)^{1/2}.\eeas
Because of the symmetry in the expression, similar inequalities hold for the sums $T_2$ and $T_4$. Therefore, 
\beas S&\leq&(2+|b_1c_1|)r^2+4(\sqrt{2}+|b_1|+|c_1|)r^3+\frac{2r^2\sqrt{-8r^{10}+31r^8-44r^6+27r^4}}{(1-r^2)^{2}}\\[2mm]
&&+\left(|b_1|+|c_1|\right)\left(\frac{r^8\left(15-20 r^2+15 r^4-4 r^6\right)}{(1-r^2)^4}-r^4\left(\log(1-r^2)+r^2\right)\right)^{1/2}\\[2mm]
&&+\frac{r^4\left(4r^2-11r+9\right)}{(1-r)^3}.\eeas
Since $|b_1|\leq2$ and $|c_1|\leq2$ and thus, we have 
\beas S&\leq &6r^2+4(\sqrt{2}+4)r^3+\frac{2r^2\sqrt{-8r^{10}+31r^8-44r^6+27r^4}}{(1-r^2)^{2}}+\frac{r^4\left(4r^2-11r+9\right)}{(1-r)^3}\\[2mm]
&&+4\left(\frac{r^8\left(15-20 r^2+15 r^4-4 r^6\right)}{(1-r^2)^4}-r^4\left(\log(1-r^2)+r^2\right)\right)^{1/2}\leq 1\eeas
for $0 <r\leq r_0$, where $r_0(\approx 0.260985)$ is the smallest positive root of the equation 
\bs\beas A(r):&=&6r^2+4(\sqrt{2}+4)r^3+\frac{2r^2\sqrt{-8r^{10}+31r^8-44r^6+27r^4}}{(1-r^2)^{2}}+\frac{r^4\left(4r^2-11r+9\right)}{(1-r)^3}\\[2mm]
&&+4\left(\frac{r^8\left(15-20 r^2+15 r^4-4 r^6\right)}{(1-r^2)^4}-r^4\left(\log(1-r^2)+r^2\right)\right)^{1/2}-1=0\eeas \es
in $(0,1)$, and it is illustrated in Figure \ref{fig5}. 
This completes the proof.\end{proof}
\begin{cor}
Let $g,h\in\mathcal{S}$ with $g''(0)=0$. Then, the function $F$ defined by (\ref{rm}) belongs to the class $\mathcal{M}$ in the disk $|z|\leq r_0$, where $r_0\thickapprox 
0.313967$ is the unique root of the equation
\bs\beas B(r):&=&2r^2+4(\sqrt{2}+2)r^3+\frac{2r^2\sqrt{-8r^{10}+31r^8-44r^6+27r^4}}{(1-r^2)^{2}}+\frac{r^4(4r^2-11r+9)}{(1-r)^3}\\
&&+2\left(\frac{r^8(15-20 r^2+15 r^4-4 r^6)}{(1-r^2)^4}-r^4\left(\log(1-r^2)+r^2\right)\right)^{1/2}-1=0\eeas\es
in $(0,1)$, as illustrated in Figure \ref{fig6}.
\begin{proof} 
The result is an immediate consequence of \textrm{Theorem \ref{T4}}, with the additional condition $b_1=0$.\end{proof}
\end{cor}
\begin{cor}
Let $g,h\in\mathcal{S}$ with $g''(0)=0$ and $h''(0)=0$. Then, the function $F$ defined by (\ref{rm}) belongs to the class $\mathcal{M}$ in the disk $|z|\leq r_0$, where $r_0\thickapprox 
0.352049$ is the unique root of the equation
\beas C(r):=2r^2+4\sqrt{2}r^3+\frac{2r^2\sqrt{-8r^{10}+31r^8-44r^6+27r^4}}{(1-r^2)^{2}}+\frac{r^4(4r^2-11r+9)}{(1-r)^3}-1=0\eeas
in $(0,1)$, as illustrated in Figure \ref{fig7}.\end{cor}
\begin{proof} 
The result is an immediate consequence of \textrm{Theorem \ref{T4}}, with the additional conditions $b_1=c_1=0$.\end{proof}
In \cite{10a}, Obradovi\'{c} {\it et al.} dealt with the problem of determining the radius of univalence of the function $G(z)$, which is defined by the quotient 
$G(z) = zf(z)/g(z)$. In this context, the functions $f(z)$ and $g(z)$ are chosen from suitable subclasses of $\mathcal{A}$. 
Our objective is to identify the largest disks for which the functions $G$, defined as the quotient of analytic functions, belong to the class $\mathcal{M}$.
\begin{theo}\label{Th8} Let $G(z)=z^2/f(z)$ for $z\in\Bbb{D}$, where $f\in\mathcal{A}$ . Then, we have the following:
\begin{itemize} 
\item[(i)] If $f\in\mathcal{S}$, then $G\in\mathcal{M}$ in the disk $|z|\leq 2-\sqrt{3}$. The result is sharp.
\item[(ii)] If $f\in\mathcal{R}(1/2)$, then $G\in\mathcal{M}$ in the disk $|z|\leq r_0$, where $r_0(\approx 0.396608)$ is the unique positive root of the equation $1-3 r+2 r^2-2 r^3\geq 0$ in $(0,1)$. The result is sharp.
\item[(iii)] If $f\in\mathcal{C}(-1/2)$, then $G\in\mathcal{M}$ in the disk $|z|\leq r_1$, where $r_1(\approx 0.304725)$ is the unique positive root of the equation $2r^4-6r^3+4r^2-4 r+1=0$ in $(0,1)$. The number $r_1$ is the best possible.
\item[(iv)] If $f\in\mathcal{G}$, then $G\in\mathcal{M}$ in the disk $|z|\leq r_2$, where $r_2(\approx 0.75085)$ is the positive root of the equation $3r-2r^2+(r^2- 5 r+4)\log(1-r)=0$ in $(0,1)$. 
\end{itemize}
\end{theo}
\begin{proof} 
\item[(i)] Let $f\in\mathcal{S}$ be such that $f(z)=z+\sum_{n=2}^\infty a_n z^n$ for $z\in\Bbb{D}$.  Then, we have $|a_n|\leq n$ for $n\geq 2$ (see \cite{r8,r9}). Since $f(z)/z\not=0$ in $\Bbb{D}$, therefore $G(z)/z\not=0$. Now 
\beas z^2\left(\frac{z}{G(z)}\right)''+G'(z)\left(\frac{z}{G(z)}\right)^2-1&=&z^3 \left(\frac{1}{G(z)}-\frac{1}{z}\right)''+z^2 \left(\frac{1}{G(z)}-\frac{1}{z}\right)'\\[2mm]
&=&z^3 \left(\frac{f(z)}{z^2}-\frac{1}{z}\right)''+z^2 \left(\frac{f(z)}{z^2}-\frac{1}{z}\right)'.\eeas
It is evident that
\beas z^2\left(\frac{f(z)}{z^2}-\frac{1}{z}\right)'=z^2\left(\frac{1}{z}+\sum_{n=2}^\infty a_n z^{n-2}-\frac{1}{z}\right)'=\sum_{n=2}^\infty (n-2)a_n z^{n-1}\quad\text{and}\\[2mm]
z^3\left(\frac{f(z)}{z}-\frac{1}{z}\right)''=z^3\left(\sum_{n=2}^\infty (n-2)a_n z^{n-3}\right)'=\sum_{n=2}^\infty (n-2)(n-3)a_n z^{n-1}.\eeas
Therefore, we have 
\bea\label{y2} \left|z^2\left(\frac{z}{G(z)}\right)''+G'(z)\left(\frac{z}{G(z)}\right)^2-1\right|\leq \sum_{n=2}^\infty (n-2)^2 |a_n| r^{n-1}\quad\text{for}\quad |z|=r.\eea
As $|a_n|\leq n$ for $n\geq 2$, therefore, we have $G\in\mathcal{M}$ if $r^2\left(3+4r-r^2\right)/(1-r)^4\leq 1$,  {\it i.e.,} if $r\leq 2-\sqrt{3}$. \\[2mm]
\indent To prove the sharpness of the radius $2-\sqrt{3}$, we consider the function $f(z)=z/(1-z)^2$. Then the quotient function $G(z)=z(1-z)^2$ and 
\beas \left|z^2\left(\frac{z}{G (z)}\right)''+G'(z)\left(\frac{z}{G(z)}\right)^2-1\right|_{z=r}=\frac{3r^2+4r^3-r^4}{(1-r)^4}>1\quad\text{for}\quad r>2-\sqrt{3},\eeas
which shows that the number $2-\sqrt{3}$ is the best possible.\\
\item[(ii)] Using similar argument as in (i), we have (\ref{y2}). Since $f\in\mathcal{R}(1/2)$, so $|a_n|\leq 1$ for $n\geq 2$. Therefore, $G\in\mathcal{M}$ if  $r^2 (1 + r)/(1-r)^3\leq 1$, {\it i.e.,} if $r\leq r_0$, where $r_0(\approx 0.396608)$ is the unique positive root of the equation $1-3 r+2 r^2-2 r^3\geq 0$ in $(0,1)$.\\[2mm]
\indent To prove the sharpness of the radius $r_0$, we consider the function $f(z)=z/(1-z)$. Therefore, $G(z)=z(1-z)$ and 
\beas \left|z^2\left(\frac{z}{G(z)}\right)''+G'(z)\left(\frac{z}{G(z)}\right)^2-1\right|_{z=r}=\frac{r^2(1+r)}{(1-r)^3}>1\quad\text{for}\quad r>r_0,\eeas
which shows that the number $r_0$ is the best possible.\\
\item[(iii)] Since $f\in\mathcal{C}(-1/2)$, so $|a_n|\leq (n+1)/2$ for $n\geq 2$. In view of (\ref{y2}), we have $G\in\mathcal{M}$ if 
$r^2 (2+2 r-r^2)/(1-r)^4\leq 1$, {\it i.e.,} $r\leq r_1$, where $r_1(\approx 0.304725)$ is the unique positive root of the equation $2r^4-6r^3+4r^2-4 r+1=0$.\\[2mm]
\indent To prove the sharpness of the radius $r_0$, we consider the function $f(z)=z(1-z/2)/(1-z)^2$. Then, we have $G(z)=z(1-z)^2/(1-z/2)$ and 
\beas \left|z^2\left(\frac{z}{G(z)}\right)''+G'(z)\left(\frac{z}{G(z)}\right)^2-1\right|_{z=r}=\frac{2r^2+2r^3-r^4}{(1-r)^4}>1\quad\text{for}\quad r>r_1,\eeas
which shows that the number $r_1$ is the best possible.\\
\item[(iv)] Since $f\in\mathcal{G}$, so $|a_n|\leq 1/(n(n-1))$ for $n\geq 2$ (see \cite{r4}). In view of (\ref{y2}), we have $G\in\mathcal{M}$ if 
$(4-3r)/(1-r)+(4-5 r+r^2)\log(1-r)/(r(1-r))\leq 1$, {\it i.e.,} $r\leq r_2$, where $r_2(\approx 0.75085)$ is the positive root of the equation $3r-2r^2+(r^2- 5 r+4)\log(1-r)=0$ in $(0,1)$, as illustrated in Figure \ref{fig8}.
\end{proof}

\begin{theo}\label{Th9} Let $H(z)=z^2/\int_0^z (t/f(t))dt$ for $z\in\Bbb{D}$, where $f\in\mathcal{A}$. Then, we have the following:
\begin{itemize} 
\item[(i)] If $f\in\mathcal{S}$, then $H\in\mathcal{M}$ in the disk $|z|\leq r_3$, where $r_3(\approx 0.7829)$ is the unique positive root of the equation 
\beas \frac{-4+4r^2+r^4}{(1-r^2)^2}- \frac{12}{r^2}\log(1-r^2)-\frac{8}{r^2}\text{Li}_2(r^2)-1=0\quad\text{in}\quad (0,1).\eeas
\item[(ii)] If $f\in\mathcal{M}$, then $H\in\mathcal{M}$ in the unit disk $\Bbb{D}$. 
\end{itemize}
\end{theo}
\begin{proof}
\item[(i)] Let $f\in\mathcal{S}$. Then, $f$ have the power series representation of the form $z/f(z)=1+\sum_{n=1}^\infty b_n z^n$. In view of \textrm{Lemma \ref{lem0}}, we have $\sum_{n=2}^\infty (n-1)|b_n|^2\leq 1$. Let $g(z)=\int_0^z (t/f(t))dt=z+\sum_{n=1}^\infty (b_n/(n+1)) z^{n+1}=z+\sum_{n=2}^\infty (b_{n-1}/n) z^{n}$. Thus, 
\beas \frac{1}{H(z)}-\frac{1}{z}=\sum_{n=2}^\infty \frac{b_{n-1}}{n} z^{n-2}=\sum_{n=1}^\infty \frac{b_{n}}{n+1} z^{n-1}.\eeas
Therefore, 
\bea\label{y3}z^2\left(\frac{z}{H(z)}\right)''+H'(z)\left(\frac{z}{H(z)}\right)^2-1&=&z^3 \left(\frac{1}{H(z)}-\frac{1}{z}\right)''+z^2 \left(\frac{1}{H(z)}-\frac{1}{z}\right)'\nonumber\\[2mm]
&=&\sum_{n=2}^\infty \frac{(n-1)^2b_{n}}{n+1} z^{n}.\eea
In view of the Cauchy-Schwarz inequality and for $|z|=r$, we have
\beas \left|z^2\left(\frac{z}{H(z)}\right)''+H'(z)\left(\frac{z}{H(z)}\right)^2-1\right|&\leq& \sum_{n=2}^\infty \frac{(n-1)^2|b_n|}{n+1} r^{n}\\
&\leq&\left(\sum_{n=2}^\infty (n-1)|b_n|^2\right)^{1/2}\left(\sum_{n=2}^\infty \frac{(n-1)^{3}}{(n+1)^2}r^{2n}\right)^{1/2}\\
&\leq& \left(\sum_{n=2}^\infty \frac{(n-1)^{3}}{(n+1)^2}r^{2n}\right)^{1/2}.\eeas
It is evident that 
\beas \sum_{n=2}^\infty \frac{(n-1)^{3}}{(n+1)^2}r^{2n}&=&\sum_{n=2}^\infty \left(-5+n+\frac{12}{(1+n)}-\frac{8}{(1+n)^2}\right)r^{2n}\\[2mm]
&=&-\frac{5r^4}{1-r^2}+\frac{r^4(2-r^2)}{(1-r^2)^2}-6 (r^2+2) - \frac{12}{r^2}\log(1-r^2)\\[2mm]
&&+2(4+r^2)-\frac{8}{r^2}\text{Li}_2(r^2)=\frac{1}{r^2} A_1(r),\eeas
where $\text{Li}_2(r^2)$ is the dilogarithm function such that 
\beas \text{Li}_2(r^2)&=&\sum_{n=1}^\infty \frac{r^{2n}}{n^2}=-\int_{0}^{r^2}\frac{\log(1-u)}{u}du\quad\text{and}\\[2mm]
A_1(r)&=&\frac{-4r^2+4r^4+r^6}{(1-r^2)^2}- 12\log(1-r^2)-8\;\text{Li}_2(r^2).\eeas
We claim that $A_1(r)\geq 0$ for $r\in[0,1]$. Differentiating $A_1(r)$ with respect to $r$, we have 
\beas A_1'(r)&=&\frac{2}{r(1-r^2)^3}\left(8 r^2 - 20 r^4 + 15 r^6 - r^8+(8-24 r^2+24 r^4 - 8 r^6)\log(1-r^2)\right)\\
&=&\frac{1}{4r(1-r^2)^6}A_2(r),\eeas
where $A_2(r)=(8 r^2 - 20 r^4 + 15 r^6 - r^8)/(8-24 r^2+24 r^4 - 8 r^6)+\log(1-r^2)$.
Therefore
\beas A_2'(r)=\frac{r^9+4r^7+r^5}{4(1-r^2)^4}\geq 0,\eeas
which shows that $A_2(r)$ is a monotonically increasing function of $r\in[0,1]$ and it follows that $A_2(r)\geq A_2(0)=0$. Thus, $A_1'(r)\geq 0$, {\it i.e.,} $A_1(r)$ is a 
monotonically increasing function of $r\in[0,1]$, and hence, we have $A_1(r)\geq A_1(0)=0$ for $r\in[0,1]$. Therefore, 
\beas \left|z^2\left(\frac{z}{H(z)}\right)''+H'(z)\left(\frac{z}{H(z)}\right)^2-1\right|\leq 1\quad\text{for}\quad r\leq r_3,\eeas
where $r_3(\approx 0.7829)$ is the unique positive root of the equation 
\beas A_3(r):=\frac{-4+4r^2+r^4}{(1-r^2)^2}- \frac{12}{r^2}\log(1-r^2)-\frac{8}{r^2}\text{Li}_2(r^2)-1=0\eeas
in $(0,1)$ and it is shown in Figure \ref{fig9}.\\
\item[(ii)] Let $f\in\mathcal{M}$. Then, $f$ can be written as $z/f(z)=1+\sum_{n=1}^\infty b_n z^n$. In view of \textrm{Lemma \ref{lem2}}, we have 
$\sum_{n=2}^\infty (n-1)^4|b_n|^2\leq 1$. Using similar argument as in (i) with the Cauchy-Schwarz inequality and from (\ref{y3}), we have 
\beas \left|z^2\left(\frac{z}{H(z)}\right)''+H'(z)\left(\frac{z}{H(z)}\right)^2-1\right|&\leq& \sum_{n=2}^\infty \frac{(n-1)^2|b_n|}{n+1} r^{n}\\
&\leq&\left(\sum_{n=2}^\infty (n-1)^4|b_n|^2\right)^{1/2}\left(\sum_{n=2}^\infty \frac{r^{2n}}{(n+1)^2}\right)^{1/2}\\[2mm]
&\leq& \left(\sum_{n=2}^\infty \frac{r^{2n}}{(n+1)^2}\right)^{1/2}.\eeas
It is evident that 
\beas \sum_{n=2}^\infty \frac{r^{2n}}{(n+1)^2}\leq \sum_{n=2}^\infty \frac{1}{(n+1)^2}=\frac{\pi^2}{6}-\frac{5}{4}<1.\eeas
Therefore, $H\in\mathcal{M}$ in the unit disk $\Bbb{D}$.
\end{proof}
\begin{theo}\label{Th5} If $f,g\in\mathcal{M}(\lambda)$ for some $0<\lambda\leq 1$ , then 
\beas \frac{fg}{(1-t)f+tg}\in\mathcal{M}(\lambda)\quad\text{for}\quad 0\leq t\leq 1.\eeas
\end{theo}
\begin{proof} As $f\in\mathcal{M}(\lambda)$, we have
\beas z^2\left(\frac{z}{f(z)}\right)''+f'(z)\left(\frac{z}{f(z)}\right)^2-1=\lambda \omega_1(z),\eeas
where $\omega_1:\mathbb{D}\to\mathbb{D}$ is analytic with $\omega_1(0)=\omega_1'(0)=0$.
It is evident that
\beas z^2\left(\frac{z}{f(z)}\right)''+f'(z)\left(\frac{z}{f(z)}\right)^2-1
=z^3 \left(\frac{1}{f(z)}-\frac{1}{z}\right)''+z^2 \left(\frac{1}{f(z)}-\frac{1}{z}\right)'.\eeas
Therefore, 
\bea\label{j4} z^3 \left(\frac{1}{f(z)}-\frac{1}{z}\right)''+z^2 \left(\frac{1}{f(z)}-\frac{1}{z}\right)'=\lambda \omega_1(z).\eea
Similarly, for $g\in\mathcal{M}(\lambda)$, we have 
\bea\label{j5} z^3 \left(\frac{1}{g(z)}-\frac{1}{z}\right)''+z^2 \left(\frac{1}{g(z)}-\frac{1}{z}\right)'=\lambda \omega_2(z),\eea
where $\omega_2:\mathbb{D}\to\mathbb{D}$ is analytic with $\omega_2(0)=\omega_2'(0)=0$.
Let $0\leq t\leq 1$ and 
\beas F_1(z)=\frac{f(z)g(z)}{(1-t)f(z)+tg(z)}.\eeas  From (\ref{j4}) and (\ref{j5}), we have 
\beas 
&&z^3 \left(\frac{1}{F_1(z)}-\frac{1}{z}\right)''+z^2 \left(\frac{1}{F_1(z)}-\frac{1}{z}\right)'\\[2mm]
&=&z^3 \left(\frac{(1-t)f(z)+tg(z)}{f(z)g(z)}-\frac{1}{z}\right)''+z^2 \left(\frac{(1-t)f(z)+tg(z)}{f(z)g(z)}-\frac{1}{z}\right)'\\[2mm]
&=&(1-t)\left(z^3 \left(\frac{1}{g(z)}-\frac{1}{z}\right)''+z^2 \left(\frac{1}{g(z)}-\frac{1}{z}\right)'\right)\\[2mm]
&&+t\left( z^3 \left(\frac{1}{f(z)}-\frac{1}{z}\right)''+z^2 \left(\frac{1}{f(z)}-\frac{1}{z}\right)'\right)\\[2mm]
&=&\lambda\omega_3(z),\eeas
where $\omega_3(z)=(1-t)\omega_1(z)+t\omega_2(z)$. It is easy to see that $\omega_3:\Bbb{D}\to\Bbb{D}$ satisfying $\omega_3(0)=0$ and $\omega_3'(0)=0$. 
Consequently, we have $F_1(z)\in\mathcal{M}(\lambda)$. This completes the proof.
\end{proof}

\section{A certain subclass of the starlike functions}
In 2017, Peng and Zhong \cite{r1} introduced the class $\Omega$ which consists of functions $f\in\mathcal{A}$ satisfying 
\bea\label{j1} |z f'(z)-f(z)|<\frac{1}{2},\quad z\in\Bbb{D}.\eea
The class $\Omega$ is a subclass of the starlike function (see \cite[Theorem 3.1]{r1}). 
As $f\in\mathcal{A}$, it follows that $f$ has the Taylor series expansion of the form $f(z)=z+\sum_{n=2}^\infty a_n z^n$ for $z\in\Bbb{D}$.
Furthermore, Peng and Zhong \cite{r1} showed that (\ref{j1}) is equivalent with 
\bea\label{j6} f(z)=z+\frac{z^2}{2}\int_0^1 \omega(zt) dt,\eea
where $\omega(z)$ is analytic in $\Bbb{D}$ with $|\omega(z)|\leq 1$ for $z\in\Bbb{D}$.
We claim that the class $\Omega$ is not a subset of the class $\mathcal{M}$.
\begin{exam}Let $f$ us consider a function defined by
\beas \frac{z}{f(z)}=1+\frac{2}{3}z+\frac{1}{3}z^3,\quad z\in\Bbb{D}.\eeas
Then, $z/f(z)\not=0$ in $\Bbb{D}$ and
\beas \left|zf'(z)-f(z)\right|=\left|z^2\left(\frac{f(z)}{z}\right)'\right|=\left|z^2\left(\frac{3}{3+2z+z^3}\right)'\right|=\left|-\frac{3z^2(3z^2+2)}{(z^3+2z+3)^2}\right|<1\eeas
for $z\in\Bbb{D}$, which shows that $f\in\Omega$. Now
\beas \left|z^2\left(\frac{z}{f(z)}\right)''+f'(z)\left(\frac{z}{f(z)}\right)^2-1\right|&=&\left|z^2\left(\frac{z}{f(z)}\right)''-z\left(\frac{z}{f(z)}\right)'+\frac{z}{f(z)}-1\right|\\
&=&\frac{4|z|^3}{3}>1\quad \text{for}\quad |z|>(3/4)^{1/3}(\approx 0.90856),\eeas
which shows that $f\not\in \mathcal{M}$.
\end{exam}
\begin{exam}Let $f$ us consider a function defined by
\beas \frac{z}{f(z)}=1+\frac{1}{2}z+\frac{1}{2}z^3,\quad z\in\Bbb{D}.\eeas
Then, $z/f(z)\not=0$ in $\Bbb{D}$. Similarly, we can show that $f\in\Omega$ while $f\not\in\mathcal{M}$.
\end{exam}
Therefore, the class $\Omega$ is not a subset of the class $\mathcal{M}$.
\section{The Bohr's phenomenon for analytic functions}
Let $f$ be an analytic function on the open unit disk $\mathbb{D}:=\{z\in\mathbb{C}:|z|<1\}$ such that $|f(z)|\leq 1$ in $\mathbb{D}$ with the series expansion  
$f(z)=\sum_{n=0}^{\infty} a_nz^n$. 
Then,
\bea\label{e2} \sum_{n=0}^{\infty}|a_n|r^n\leq 1\quad\text{for}\quad|z|=r\leq\frac{1}{3}.\eea
Note that if $|f(z)|\leq 1$ in $\mathbb{D}$ and $|f(a)|=1$ for some $a\in\mathbb{D}$, then $f(z)$ reduces to a unimodular constant function (see \cite[Strict Maximum Principle (Complex Version), P. 88]{3a}).
In this context, the quantity $1/3$ is known as the Bohr radius and it can't be improved. The inequality (\ref{e2}) is known as the Bohr inequality.
In fact, H. Bohr \cite{5t} was the first to derive the inequality (\ref{e2}) for $r\leq 1/6$. However, this result was subsequently refined by Weiner, Riesz, and Schur \cite{200} to yield the improved value of $1/3$.
The Bohr inequality (\ref{e2}) can be written as 
\beas\sum_{n=1}^{\infty}|a_n|r^n\leq 1-|a_0|=1-|f(0)|=d\left(f(0),\pa\Bbb{D}\right)\quad\text{for}\quad|z|=r\leq1/3,\eeas
where $d\left(f(0),\pa\Bbb{D}\right)$ denotes the Euclidean distance between $f(0)$ and the boundary $\pa\Bbb{D}$ of $\Bbb{D}$.
The concept of Bohr's phenomenon can be generalized to the class $\mathcal{F}$ 
consisting of analytic functions $f$ that map from $\mathbb{D}$ into a given domain $U\subseteq\mathbb{C}$ such that $f(\mathbb{D})\subseteq U$. The 
generalized Bohr phenomenon for the class $\mathcal{F}$ is defined as:\\[1mm]
For $U\subset \Bbb{C}$, we have to find the largest number $r_U\in(0,1)$ such that 
\beas  \sum_{n=1}^{\infty}|a_n|r^n\leq d\left(f(0),\pa f(\Bbb{D})\right)\quad\text{holds for}\quad|z|=r\leq r_U.\eeas
We refer to \cite{B1,B2,B3,B4,B5,B6,B7,B8,BM2024} and the references listed therein for an in-depth investigation on several other aspects of Bohr's inequality. Besides the concept of Bohr 
radius, there is another concept called Rogosinski radius (see \cite{24t}). Motivated by the Rogosinski radius, Kayumov and Ponnusamy \cite{300} have considered the Bohr-Rogosinski sum $R_N^f(z)$ defined by
\beas R_N^f(z):=|f(z)|+\sum_{n=N}^{\infty}|a_n||z|^n.\eeas
The Bohr-Rogosinski sum
$R_N^f(z)$ is related to the classical Bohr sum (Majorant series) in which $N=1$ and $f(z)$ is replaced by $f(0)$. For an analytic function $f$ in $\Bbb{D}$ with $|f(z)|<1$
in $\Bbb{D}$, Kayumov and Ponnusamy \cite{300} have defined the Bohr-Rogosinski radius as the largest number $r\in (0, 1)$ such that $R_N^f(z)\leq 1$ holds for $|z|<r$.\\[2mm]
\indent 
The following lemmas pertain to the growth estimate for functions in the class $\Omega$.
\begin{lem}\cite{r1}\label{lem8} If $f\in\Omega$, then
\beas |z|-\frac{1}{2}|z|^2\leq |f(z)|\leq |z|+\frac{1}{2}|z|^2\quad\text{and}\quad 1-|z|\leq |f'(z)|\leq 1+|z|.\eeas
For each $z\in\Bbb{D}$, $z\not=0$, equality occurs in both estimates if and only if $f(z)=z+\frac{1}{2}\lambda z^2$, where $|\lambda|=1$. 
\end{lem}
The following result provides a sufficient condition for functions to be members of the class $\Omega$.
\begin{lem}\label{lem9} Let $f\in\mathcal{A}$ be given by $f(z)=z+\sum_{n=2}^\infty a_n z^n$ and $\sum_{n=2}^\infty (n-1)|a_n|\leq 1/2$. Then, $f\in\Omega$. \end{lem}
\begin{proof}
As $f(z)=z+\sum_{n=2}^\infty a_n z^n$ for $z\in\Bbb{D}$, thus
\beas \left|zf'(z)-f(z)\right|=\left|\sum_{n=2}^\infty (n-1)a_n z^n \right|\leq\sum_{n=2}^\infty (n-1)|a_n||z|^n<\sum_{n=2}^\infty (n-1)|a_n|\leq\frac{1}{2}. \eeas
Therefore, $f\in\Omega$ and this completes the proof.
\end{proof}
We define the class 
\beas \Omega_\mathcal{A}:=\left\{f(z)=z+\sum_{n=2}^\infty a_n z^n\in\mathcal{A}: \sum_{n=2}^\infty (n-1)|a_n|\leq 1/2\right\}.\eeas
In view of \textrm{Lemma \ref{lem9}}, we have $\Omega_\mathcal{A}\subseteq\Omega$.
\section{The Bohr phenomenon for the class $\Omega_\mathcal{A}$}
In the following result, we obtain the sharp Bohr-Rogosinski radius for functions in the class $\Omega_\mathcal{A}$.
\begin{theo}
Let $f\in\Omega_\mathcal{A}$ be given by $f(z)=z+\sum_{n=2}^\infty a_n z^n$ for $z\in\Bbb{D}$. Then, we have 
\beas |f(z)|+\sum_{n=2}^\infty |a_n| r^n\leq d\left(f(0),\pa f(\Bbb{D})\right)\quad\text{for}\quad r\leq (\sqrt{3}-1)/2,\eeas
The radius $(\sqrt{3}-1)/2$ is the best possible.
\end{theo}
\begin{proof}
Let $f\in\Omega_\mathcal{A}$ be such that $f(z)=z+\sum_{n=2}^\infty a_n z^n$ for $z\in\Bbb{D}$. Using \textrm{Lemma \ref{lem8}}, we have
\bea\label{t1} |z|-\frac{1}{2}|z|^2\leq |f(z)|, \quad\text{where}\quad|z|<1.\eea
The Euclidean distance between $f(0)$ and $\partial f(\mathbb{D})$, the boundary of $f(\mathbb{D})$ is $d\left(f(0),\partial f(\mathbb{D})\right):=\liminf_{|z|\to 1^-}|f(z)-f(0)|=\liminf_{|z|\to 1^-}|f(z)|$, as $f(0)=0$. 
From (\ref{t1}), we get
\bea\label{t2}\frac{1}{2}\leq d\left(f(0),\partial f(\mathbb{D})\right).\eea
From the definition of the class $\Omega_\mathcal{A}$, we have 
\bea\label{t3} &&\sum_{n=2}^\infty|a_n|\leq \sum_{n=2}^\infty (n-1)|a_n|\leq 1/2\nonumber\\[2mm]\text{and}
&&\sum_{n=2}^{\infty}|a_n|r^n=r^2\left(|a_2|+|a_3|r+|a_4|r^2+\cdots\right)\leq r^2\sum_{n=2}^\infty|a_n|\leq \frac{r^2}{2}.\eea
In view of \textrm{Lemma \ref{lem8}} and by the use of inequality (\ref{t3}), we have
\beas|f(z)|+\sum_{n=2}^{\infty}|a_n|r^n\leq r+r^2\leq \frac{1}{2}\quad\text{for}\quad 0<r\leq r_0, \eeas
where $r_0\in(0,1)$ is the unique root of the equation $G_1(r)=0$, where $G_1:[0,1)\to\mathbb{R}$ is a function defined by
\beas G_1(r)=2r+2r^2-1.\eeas
It is clear that $G_1$ is continuous in $[0,1)$ and strictly increasing function of $r\in[0,1)$ with $G_1(0)=-1<0$ and $\lim_{r\to 1^-}G_1(r)=3$. 
By the Intermediate Value Theorem, the equation $G_1(r)=0$ has a unique root $r_0=(\sqrt{3}-1)/2$ in $(0,1)$.
From (\ref{t2}), for $|z|=r\leq r_0=(\sqrt{3}-1)/2$, we have
\beas |f(z)|+\sum_{n=2}^{\infty}|a_n|r^n\leq d\left(f(0),\partial f(\mathbb{D})\right).\eeas
In order to prove the sharpness of the constant $r_0$, we consider the function
\bea\label{t6} f_1(z)=z+z^2/2.\eea
Note that $f_1(0)=0$ and $f_1\in\Omega_\mathcal{A}$. For $z=r$, we have 
\bea\label{t7}|f_1(r)-f_1(0)|=|f_1(r)|=\left| r+\frac{r^2}{2}\right|=r+\frac{r^2}{2}\quad\text{and}\quad\liminf_{r\to 1^{-}}|f_1(r)|=\frac{3}{2}.\eea
and for $z=-r$, we have 
\bea\label{t8}|f_1(-r)-f_1(0)|=\left|-r+\frac{(-r)^2}{2}\right|=r-\frac{r^2}{2}\quad\text{and}\quad\liminf_{r\to 1^{-}}|f_1(-r)|=\frac{1}{2}.\eea 
From (\ref{t7}) and (\ref{t8}), we have 
\bea\label{t9}d\left(f_1(0),\partial f_1(\mathbb{D})\right)=\frac{1}{2}.\eea
For $f=f_1$ and $z=r$, a simple computation using (\ref{t9}) shows that
\beas |f_1(r)|+\sum_{n=2}^{\infty}|a_n|r^n=r+r^2>\frac{1}{2}=d\left(f_1(0),\partial f_1(\Bbb{D})\right)\quad\text{for}\quad r>r_0,\eeas
which shows that the radius $r_0$ is the best possible. This completes the proof.
\end{proof}
In the following result, we obtain the sharp Bohr radius for functions in the class $\Omega_\mathcal{A}$.
\begin{theo}
Let $f\in\Omega$ be given by $f(z)=z+\sum_{n=2}^\infty a_n z^n$ for $z\in\Bbb{D}$. Then, we have 
\beas r+\sum_{n=2}^\infty |a_n| r^n\leq d\left(f(0),\pa f(\Bbb{D})\right)\quad\text{for}\quad r\leq \sqrt{2}-1.\eeas
The radius $\sqrt{2}-1$ is the best possible.
\end{theo}
\begin{proof}
Let $f\in\Omega_\mathcal{A}$ be such that $f(z)=z+\sum_{n=2}^\infty a_n z^n$ for $z\in\Bbb{D}$. Using similar argument as in the proof of \textrm{Theorem \ref{Th5}} and in view of \textrm{Lemma \ref{lem8}}, we have the inequality (\ref{t2}) and (\ref{t3}). 
In view of inequality (\ref{t3}), we have
\beas r+\sum_{n=2}^{\infty}|a_n|r^n\leq r+\frac{r^2}{2}\leq \frac{1}{2}\eeas
for $0<r\leq r_1=\sqrt{2}-1$, where $r_1\in(0,1)$ is the unique root of the equation $G_2(r):=2r+r^2-1=0$.
From (\ref{t2}), for $|z|=r\leq r_1=\sqrt{2}-1$, we have
\beas r+\sum_{n=2}^{\infty}|a_n|r^n\leq d\left(f(0),\partial f(\mathbb{D})\right).\eeas
In order to prove the sharpness of the constant $r_1$, we consider the function $f_1$ defined in (\ref{t6}).
It is evident that $f_1\in\Omega_\mathcal{A}$ with $f_1(0)=0$ and $d\left(f_1(0),\partial f_1(\mathbb{D})\right)=1/2$.
For $|z|=r_1$, we have $r+\sum_{n=2}^{\infty}|a_n|r^n=r_1+r_1^2/2=1/2=d\left(f_1(0),\partial f_1(\Bbb{D})\right)$.
Thus, the radius $r_1$ is the best possible. This completes the proof. \end{proof}
In the following result, we obtain the sharply improved version of Bohr inequality for functions in the class $\Omega_\mathcal{A}$.
\begin{theo}
Let $f\in\Omega$ be given by $f(z)=z+\sum_{n=2}^\infty a_n z^n$ for $z\in\Bbb{D}$. Then, we have 
\beas |f(z)|+|f'(z)||z|+\sum_{n=2}^\infty |a_n| r^n\leq d\left(f(0),\pa f(\Bbb{D})\right)\quad\text{for}\quad r\leq (\sqrt{2}-1)/2,\eeas
The number $(\sqrt{2}-1)/2$ is the best possible. 
\end{theo}
\begin{proof}
Let $f\in\Omega_\mathcal{A}$ be such that $f(z)=z+\sum_{n=2}^\infty a_n z^n$ for $z\in\Bbb{D}$. Using similar argument as in the proof of \textrm{Theorem \ref{Th5}}and in view of \textrm{Lemma \ref{lem8}}, we have the inequality (\ref{t2}) and (\ref{t3}). 
In view of \textrm{Lemma \ref{lem8}} and by the use of inequality (\ref{t3}), we have
\beas|f(z)|+|f'(z)||z|+\sum_{n=2}^{\infty}|a_n|r^n\leq 2r+2r^2\leq \frac{1}{2}\quad\text{for}\quad 0<r\leq r_2, \eeas
where $r_2\in(0,1)$ is the unique root of the equation $G_3(r):=4r+4r^2-1=0$.
From (\ref{t2}), for $|z|=r\leq r_2=(\sqrt{2}-1)/2$, we have
\beas |f(z)|+|f'(z)||z|+\sum_{n=2}^{\infty}|a_n|r^n\leq d\left(f(0),\partial f(\mathbb{D})\right).\eeas
In order to prove the sharpness of the constant $r_2$, we consider the function $f_1$ defined in (\ref{t6}).
It is evident that $f_1\in\Omega_\mathcal{A}$ with $f_1(0)=0$ and $d\left(f_1(0),\partial f_1(\mathbb{D})\right)=1/2$.
For $f=f_1$ and $z=r$, we have 
\beas |f_1(r)|+|f_1'(r)|r+\sum_{n=2}^{\infty}|a_n|r^n=2r+2r^2>\frac{1}{2}=d\left(f_1(0),\partial f_1(\Bbb{D})\right)\quad\text{for}\quad r>r_2,\eeas
which shows that the radius $r_2$ is the best possible. This completes the proof. \end{proof}
\section{Declarations}
{\bf Acknowledgment:} The work of the first author is supported by University Grants Commission (IN) fellowship (No. F. 44 - 1/2018 (SA - III)).\\
{\bf Conflict of Interest:} The authors declare that there are no conflicts of interest regarding the publication of this paper.\\
{\bf Availability of data and materials:} Not applicable.


\begin{thebibliography}{33}
\bibitem{1} {\sc L. A. Aksent\'ev}, Sufficient conditions for univalence of regular functions (Russian), {\it Izv. Vys\v{s}. Uc\v{e}bn. Zaved. Matematika} {\bf3} (1958), 3--7.
\bibitem{1a} {\sc L. A. Aksent\'ev} and  {\sc F. G. Avhadiev}, A certain class of univalent functions (Russian), {\it Izv. Vys\v{s}. Uc\v{e}bn. Zaved. Matematika} {\bf10} (1970), 12--20.
\bibitem{B3}{\sc S. A. Alkhaleefah, I. R. Kayumov} and  {\sc S. Ponnusamy}, On the Bohr inequality with a fixed zero coefficient, {\it Proc. Am. Math. Soc.} {\bf147}(12) (2019), 5263--5274.
\bibitem{B4}{\sc S. A. Alkhaleefah, I. R. Kayumov} and  {\sc S. Ponnusamy}, Bohr-Rogosinski inequalities for bounded analytic functions, {\it Lobachevskii J. Math.} {\bf41} (2020), 2110--2119.
\bibitem{B5}{\sc M. B. Ahamed, V. Allu} and  {\sc H. Halder}, The Bohr phenomenon for analytic functions on a shifted disk, {\it Ann. Fenn. Math.} {\bf47} (2022), 103--120.
\bibitem{V1}{\sc  V. Allu, V. Arora} and  {\sc A. Shaji}, On the Second Hankel Determinant of Logarithmic Coefficients for Certain Univalent Functions, {\it Mediterr. J. Math.} {\bf20} (2023), 81.
\bibitem{B6}{\sc V. Allu} and  {\sc H. Halder}, Bohr phenomenon for certain close-to-convex analytic functions, {\it Comput. Methods Funct. Theory} {\bf22} (2022), 491--517
\bibitem{r8}{\sc L. De Branges}, A proof of the Bieberbach conjecture, {\it Acta Math.} {\bf 154} (1985), 137--152.
\bibitem{5t} {\sc H. Bohr}, A theorem concerning power series, {\it Proc. London Math. Soc.} {\bf13}(2) (1914), 1--5.
\bibitem{103} {\sc R. Biswas}, Second Hankel determinant of logarithmic coefficients for $\mathcal{G}(\alpha)$ and $\mathcal{P}(M)$, {\it J. Anal.} 32(5) (2024), 3019--3037.
\bibitem{BM2024}{\sc R. Biswas} and {R. Mandal}, The Bohr's Phenomenon for the class of $K$-quasiconformal harmonic mappings, preprint, https://doi.org/10.48550/arXiv.2411.03352.
\bibitem{200} {\sc P. G. Dixon}, Banach algebras satisfying the non-unital von Neumann inequality, {\it Bull. Lond. Math. Soc.} {\bf27}(4) (1995), 359--362.
\bibitem{2}{\sc P. L. Duren}, Univalent Functions, {\it Grundlehren der mathematischen Wissenschaften}, vol. {\bf259}, Springer-Verlag, New York, Berlin, Heidelberg, Tokyo, 1983.
\bibitem{B7}{\sc S. Evdoridis, S. Ponnusamy} and  {\sc A. Rasila}, Improved Bohr's inequality for shifted disks, {\it Results Math} {\bf76} (2021), 14.
\bibitem{20} {\sc B. Friedman}, Two theorems on schlicht functions, {\it Duke Math. J.} {\bf13}(2) (1946), 171--177.
\bibitem{3} {\sc R. Fournier} and  {\sc S. Ponnusamy}, A class of locally univalent functions defined by a differential inequality, {\it Complex Var. Elliptic Equ.} {\bf52}(1) (2007), 1--8.
\bibitem{3a}{\sc T. W. Gamelin}, Complex Analysis, Springer-Verlag, New York, 2000.
\bibitem{4} {\sc A. W. Goodman}, Univalent Functions, Vols. {\bf1-2}, Mariner, Tampa, Florida, 1983.
\bibitem{r7} {\sc I. Jovanvi\'c} and  {\sc M. Obradovi\'c}, A note on certain classes of univalent functions, {\it Filomat} {\bf9}(1) (1995), 69--72.
\bibitem{300} {\sc I. R. Kayumov} and  {\sc S. Ponnusamy}, Bohr-Rogosinski radius for analytic functions, preprint, https://doi.org/10.48550/arXiv.1708.05585.
\bibitem{B8}{\sc G. Liu, Z. H. Liu} and  {\sc S. Ponnusamy}, Refined Bohr inequality for bounded analytic functions, {\it Bull. Sci. Math.} {\bf173} (2021), 103054.
\bibitem{B1} {\sc R. Mandal, R. Biswas }and  {\sc S. K. Guin}, Geometric studies and the Bohr radius for certain normalized harmonic mappings, {\it Bull. Malays. Math. Sci. Soc.} {\bf47} (2024), 131.
\bibitem{15} {M. Nunokawa, S. Ponnusamy} and  {\sc S. Owa}, One criterion for univalency, {\it Proc. Amer. Math. Soc.} {\bf106} (1989), 1035--1037.
\bibitem{5} {\sc M. Obradovi\'{c}} and  {\sc S. Ponnusamy}, New criteria and distortion theorems for univalent functions, {\it Complex Var. Theory Appl.} {\bf44} (2001), 173--191.
\bibitem{19} {\sc M. Obradovi\'{c}} and  {\sc S. Ponnusamy}, Radius properties for subclasses of univalent functions, {\it Analysis} {\bf25} (2005), 183--188.
\bibitem{6} {\sc M. Obradovi\'{c}} and  {\sc S. Ponnusamy}, Univalence and starlikeness of certain integral transforms defined by convolution of analytic functions, {\it J. Math. Anal. Appl.} {\bf336} (2007), 758--767.
\bibitem{9} {\sc M. Obradovi\'{c}} and  {\sc S. Ponnusamy}, On certain subclasses of univalent functions and radius properties, {\it Rev. Roumaine Math. Pures Appl.} {\bf54}(4) (2009), 317-329.
\bibitem{11} {\sc M. Obradovi\'{c}} and  {\sc S. Ponnusamy}, Coefficient characterization for certain classes of univalent functions, {\it Bull. Belg. Math. Soc. Simon Stevin} {\bf16} (2009), 251--263.
\bibitem{13} {\sc M. Obradovi\'{c}} and  {\sc S. Ponnusamy}, A class of univalent functions defined by a differential inequality, {\it Kodai Math. J.} {\bf30} (2011), 169--178.
\bibitem{7} {\sc M. Obradovi\'{c}} and  {\sc S. Ponnusamy}, On a class of univalent functions, {\it Appl. Math. Lett.} {\bf25} (2012), 1373--1378.
\bibitem{8} {\sc M. Obradovi\'{c}} and  {\sc S. Ponnusamy}, Product of univalent functions, {\it Math. Comput. Model.} {\bf57} (2013), 793--799.
\bibitem{17} {\sc M. Obradovi\'{c}} and  {\sc S. Ponnusamy}, Criteria for univalent functions in the unit disk, {\it Arch. Math.} {\bf100} (2013), 149--157.
\bibitem{10}{\sc M. Obradovi\'{c}, S. Ponnusamy, V. Singh} and  {\sc P. Vasundhra}, Univalency, starlikesess and convexity applied to certain classes of rational functions, {\it Analysis} {\bf22}(3) (2002), 225--242.
\bibitem{10a} {\sc M. Obradovi\'{c}, S. Ponnusamy} and  {\sc N. Tuneski}, Radius of univalence of certain combination of univalent and analytic functions, {\it Bull. Malays. Math. Sci. Soc.} {\bf35}(2) (2012), 325--334.
\bibitem{18} {\sc M. Obradovi\'{c}, S. Ponnusamy} and  {\sc K. J. Wirths}, Geometric Studies on the Class $\mathcal{U}(\lambda)$, {\it Bull. Malays. Math. Sci. Soc.} {\bf39} (2016), 1259--1284.
\bibitem{r4} {\sc M. Obradovi\'{c}, S. Ponnusamy} and  {\sc K. J. Wirths}, Coefficient characterizations and sections for some univalent functions, {\it Sib. Math. J.} {\bf54}(4) (2013), 1259--1284.
\bibitem{r6} {\sc S. Ozaki}, On the theory of multivalent functions II, {\it Sci. Rep. Tokyo Bunrika Daigaku. Sect. A} {\bf 4} (1941), 45--87.
\bibitem{14} {\sc S. Ozaki} and  {\sc M. Nunokawa}, The Schwarzian derivative and univalent functions, {\it Proc. Amer. Math. Soc.} {\bf33} (1972), 392--394.
\bibitem{r5} {\sc S. Ponnusamy} and  {\sc V. Allu}, Region of variability of two subclasses of univalent functions, {\it J. Math. Anal. Appl.} {\bf 332}(2) (2007), 1322--1333. 
\bibitem{r1} {\sc Z. Peng} and  {\sc G. Zhong}, Some properties for certain classes of univalent functions defined by differential inequalities, {\it Acta Math. Sci.} {\bf 37B}(1) (2017), 69--78.
\bibitem{B2}{\sc S. Ponnusamy, R. Viajayakumar} and  {\sc K. -J. Wirths}, New inequalities for the coefficients of unimodular bounded functions, {\it Results Math} {\bf75} (2020), 107.
\bibitem{12} {\sc M. O. Reade, H. Silverman} and  {\sc P. G. Todorov}, On the starlikeness and convexity of a class of analytic functions, {\it Rend. Circ. Mat. Palermo} {\bf 33} (1984), 265--272.
\bibitem{24t} {\sc W. Rogosinski}, \"Uber Bildschranken bei Potenzreihen und ihren Abschnitten, {\it Math. Z.} {\bf17} (1923), 260--276.
\bibitem{r9} {\sc P. Todorov},  A simple proof of the Bieberbach conjecture, {\it Bull. Cl. Sci., VI. S\'er., Acad. R. Belg.} {\bf3}(12) (1992), 335--356.
\bibitem{r3} {\sc D. K. Thomas, N. Tuneski} and  {\sc V. Allu}, Univalent functions : A Primer, {\it De Gruyter Studies in Mathematics}, {\bf 69}, De Gruyter, Berlin, 2018.
\end{thebibliography}
\end{document}